\documentclass[article, 11pt,reqno]{amsart}
\usepackage[english]{babel}
\usepackage[utf8]{inputenc}
\usepackage{thm-restate}
\usepackage{amsmath, amssymb, amsthm}
\usepackage[table]{xcolor}

\usepackage{graphics}
\usepackage{todonotes}
\usepackage{physics}
\usepackage{setspace}
\usepackage[margin=0.8in]{geometry}
\usepackage[hidelinks]{hyperref}

\usepackage{caption}
\usepackage{subcaption}
\usepackage{floatrow}
\usepackage{xpatch}

\usepackage{graphicx}
\usepackage{mathrsfs}
\usepackage{verbatim}

\DeclareMathOperator*{\UEG}{UEG}

\newcommand{\Z}{\mathbb{Z}}
\newcommand{\arctanh}{\text{arctanh}}
\newcommand{\cc}{\leftrightarrow}
\newcommand{\Prb}{\mathbb{P}}
\newcommand{\Prbcur}{\mathbf{P}}

\newcommand{\T}{\mathbb T}
\renewcommand{\P}{\mathbb P}
\newcommand{\Even}{\Omega_{\emptyset}}

\usepackage{tikz}
\usepackage{pgfplots}
\usepackage{tikz-3dplot} 
\usetikzlibrary{positioning}
\usetikzlibrary{arrows}

\newcommand{\M}{\mathbb{M}}
\newcommand{\id}{1\! \!1}

\newtheorem{theorem}{Theorem}[section]
\newtheorem{definition}[theorem]{Definition}
\newtheorem{lemma}[theorem]{Lemma}

\newtheorem{remark}[theorem]{Remark}
\newtheorem{corollary}[theorem]{Corollary}

\newtheorem{question}[theorem]{Question}
\newtheorem{proposition}[theorem]{Proposition}

\usepackage{cleveref}

\newlength\tindent
\pgfplotsset{compat=1.18}
\title[Non-uniqueness for graphical representations of the Ising]{Non-uniqueness of phase transitions for graphical representations of the Ising model on tree-like graphs}

 \author{Ulrik Thinggaard Hansen}
\address{Ulrik Thinggaard Hansen \\ QMATH, Department of Mathematical Sciences, University of Copenhagen, Universitetsparken 5, 2100 Copenhagen, Denmark}
\email{utha@math.ku.dk}

\author{Frederik Ravn Klausen}
\address{Frederik Ravn Klausen \\ Princeton University, Department of Mathematics, Fine Hall}
\email{fk3206@princeton.edu}

\author{Peter Wildemann}
\address{Peter Wildemann, University of Cambridge, Statistical Laboratory, DPMMS}
\email{pw477@cam.ac.uk}

\begin{document}
\begin{abstract}
We consider the graphical representations of the Ising model on tree-like graphs. We construct a class of graphs on which the loop $\mathrm{O}(1)$ model and the single random current exhibit a non-unique phase transition with respect to the inverse temperature, highlighting the non-monotonicity of both models. It follows from the construction that there exist infinite graphs $\mathbb{G}\subseteq \mathbb{G}'$ such that the uniform even subgraph of $\mathbb{G}'$ percolates and the uniform even subgraph of $\mathbb{G}$ does not.
We also show that on the wired $d$-regular tree, the phase transitions of the loop $\mathrm{O}(1)$, the single random current, and the random-cluster models are all unique and coincide. 
\end{abstract}
\maketitle

\section{Introduction} The Ising model needs no introduction as one of the cornerstones of statistical mechanics, and over the past 50 years its so-called graphical representations have become one of the main tools for its rigorous study \cite{aizenman1982geometric, DC16, DC17}. 
Consequently, they are increasingly regarded as objects of independent interest \cite{chen2023conformal,duminil2021conformal2,duminil2021conformal, Gri06, hansen2023uniform,hansen2022strict}. The most prominent of these is the \emph{random-cluster model}\footnote{We only consider the case of cluster weight $q=2$ in which case it is also referred to as the FK-Ising model.} $\phi_x$, introduced in \cite{fortuin1972random}. The \emph{loop $\mathrm{O}(1)$\! model} $\ell_x$, was introduced by Van der Waerden \cite{van1941lange} as the high-temperature expansion of the Ising model. On finite graphs,  $\ell_x$ can be defined as Bernoulli percolation $\mathbb{P}_{p}$ at parameter $p=\frac{x}{1+x}$ conditioned on being even (that is every vertex has even degree).  The \emph{random current} representation $\Prbcur_x$ was first introduced in \cite{griffiths1970concavity} and given a useful probabilistic interpretation in \cite{aizenman1982geometric}. While usually considered as a multigraph, we will be concerned with its traced version (i.e.\ the induced simple graph), which can be defined as $\Prbcur_x = \ell_x \cup \P_{1- \sqrt{1-x^2}}$, where $\mu_1 \cup \mu_2 $ is the distribution of $\omega_1\cup \omega_2$ under $\mu_1\otimes \mu_2$ ($\omega_i \sim \mu_i$).
We also refer to $\Prbcur_x \cup \Prbcur_x$ as the \emph{double random current}.

Similarly, we note that the random-cluster can be defined as $\phi_x = \ell_x \cup \P_{x}$ by \cite{evertz2002new,grimmett2007random}.
As all models are obtained from the loop $\mathrm{O}(1)$ model via "sprinkling" with Bernoulli percolation, much of our analysis will focus on the former.

In this paper, we investigate a number of natural questions regarding the interplay between the graphical representations which, while by no means ground-breaking, offer some conceptual clarification which is at present not well-represented in the literature.

\subsection{Results}
In this paper, we prove non-uniqueness of the percolative phase transition of the (free) loop $\mathrm{O}(1)$ model $\ell_{x, \mathbb{G}}^0$ (the infinite volume measure is defined in \eqref{eq:loop-cluster-coupling} below).
Here, the index $\mathbb{G}$ denotes the underlying graph and the superscript $0$ denotes free boundary conditions. We also let $\mathcal{C}_{\infty}$ denote the event that there exists an infinite cluster.
\begin{restatable}{theorem}{nonuniqueness}
     \label{thm:non_uniqueness}
There exists a graph $\mathbb{M}$ where $x \mapsto \ell^0_{x,\mathbb{M}}[\mathcal{C}_{\infty}]$ is not monotone. 
\end{restatable}

In \Cref{cor:random_current_non_monotone} we prove the same result for the (traced, sourceless) single random current model $\Prbcur_{x}^0$.
Next, consider the \emph{uniform even subgraph} $\UEG_{\mathbb{G}}$, defined to be the uniform measure on even subgraphs of $\mathbb{G}$.
This model is intimately related to the Ising model \cite{angel2021uniform, grimmett2007random,  hansen2023uniform} and can be understood as a special case of $\ell_{x, G}^0$ for $x=1$.
Using \Cref{thm:non_uniqueness}, we prove that percolation of the uniform even subgraph is not monotone in the graph. The result holds for both the free and wired infinite volume measures, both of which will be defined below. 

 \begin{restatable}{corollary}{counterexample}
 \label{cor:counterexample}
 There exist graphs $\mathbb{G}' \subset \mathbb{G}$ such that $\UEG_\mathbb{G}[\mathcal{C}_\infty] = 0$ and $\UEG_{\mathbb{G}'}[\mathcal{C}_\infty] = 1$.
\end{restatable}

In \cite{hansen2023uniform}, it was proven that on the hypercubic lattice $\Z^d$, the regime of exponential decay for the loop $O(1)$ and random current model coincides with the high temperature phase of the Ising model. Here we establish the phase diagrams for the $d$-regular tree $\mathbb{T}^d$ as well as the graph obtained from the $d$-regular tree by replacing every edge with a cycle of length $2n$ (and glued through opposite points of the cycle),  henceforth denoted $\mathtt{C}_{n}^{d}$. For a boundary condition $\xi \in \{0,1\}$, corresponding respectively to free and wired boundaries, we define the critical point of the loop $O(1)$ model via $x_c(\ell^{\xi}_{\mathbb{G}}) = \inf_{x \in [0,1]} \{\ell^{\xi}_{x,\mathbb{G}}[\mathcal{C}_\infty] >0 \}$.
Analogous definitions are used for the other models.

\begin{theorem} \label{thm:main_theorem_equal_different}
    For any $d\geq 2$ and $n\geq 1,$ it holds that
\begin{align*}
& x_c(\ell^1_{\mathbb{T}^d})  
 =  x_c(\Prbcur^1_{\mathbb{T}^d}) 
  =  x_c(\Prbcur^1_{\mathbb{T}^d}\cup \Prbcur^1_{\mathbb{T}^d}) 
 = x_c(\phi^1_{\mathbb{T}^d}) \\
& x_c(\ell^0_{\mathbb{T}^d})  
 > x_c(\Prbcur^0_{\mathbb{T}^d}) 
  >  x_c(\Prbcur^0_{\mathbb{T}^d}\cup \Prbcur^0_{\mathbb{T}^d}) 
 > x_c(\phi^0_{\mathbb{T}^d}).
 \end{align*} 
 The same statements are true for the graph $\mathtt{C}_{n}^{d}$. In both cases, all phase transitions are unique.
\end{theorem}
This theorem is the most basic illustration of the mechanism first investigated in \cite{GMM18}: While it is surprising that the single random current, double random current and random-cluster model should share a single phase transition\footnote{See the introduction in \cite{hansen2023uniform} for further explication.},  this phenomenon ultimately boils down to the existence of long loops in the ambient graph. 
In the absence of loops (as in the free tree), the phase transitions should be distinct, and if there are only long loops (as in the wired tree), we expect them to be one and the same. This carries over to the situation where all loops are of uniformly bounded length (as for the free loops on $\mathtt{C}^d_n$).

In Figure \ref{fig:Phase diagram}, we provide a graphical overview of the results. Finally, we prove that the phase transitions of the uniform even subgraph and Bernoulli percolation are not in any way related:

\begin{restatable}{theorem}{norelation} \label{norelation}
For any $\varepsilon>0,$ there exists a graph $\mathbb{G}^{\varepsilon}$ with $p_c(\mathbb{P}_{p,\mathbb{G}^{\varepsilon}})\in (1-\varepsilon,1)$ and $\UEG_{\mathbb{G}^\varepsilon}[\mathcal{C}_{\infty}]=1$.
\end{restatable}

\begin{figure}
\begin{center} 
\begin{tikzpicture}
\draw[dashed] (0,0)--(5,0);
\draw[line width=0.5mm] (5,0)--(10,0);
\node (X) at (-1.5,0) { $\ell_{\T^d}^1$, $\Prbcur_{\T^d}^1$, $\phi_{\T^d}^1$   }; 
\draw node[fill, circle, scale=0.6]  at (5,0) {};

\node (X) at (-1.5 , -2) { $\ell_{\mathtt{C}_{n}^{d}}^0$, $\Prbcur_{\mathtt{C}_{n}^{d}}^0$, $\phi_{\mathtt{C}_{n}^{d}}^0$  
 }; 
\draw[dashed] (0,-2)--(4,-2);
\draw[thin] (4,-2)--(8,-2);
\draw[line width=0.5mm] (8,-2)--(10,-2);
\node (X) at (4, -1.5) {$x_c(\phi^0)$}; 
\draw node[fill,circle, scale=0.6]  at (4,-2) {};
\node (X) at (6, -1.5) {$x_c(\Prbcur^0)$}; 
\draw node[fill,circle, scale=0.6]  at (6,-2) {};
\node (X) at (8, -1.5) {$x_c(\ell^0)$}; 
\draw node[fill,circle, scale=0.6]  at (8,-2) {};

\node (X) at (2,0.5) {exp decay}; 
\node (X) at (8,0.5) {percolation}; 
\node (X) at (2,-1.5) {exp decay}; 

\node (X) at (0,0.5) {$x=0$}; 
\node (X) at (10,0.5) {$x=1$}; 

\node (X) at (-1, -4) {$\ell_\M^0$}; 
\draw[dashed] (0,-4)--(4,-4);
\draw[dashed] (4.5,-4)--(5,-4);

\draw node[fill,circle, scale=0.6]  at (5,-4) {};
\draw node[fill,circle, scale=0.6]  at (4.5,-4) {};
\draw[line width=0.5mm] (4,-4)--(4.5,-4);
\draw node[fill,circle, scale=0.6]  at (4,-4) {};
\draw[line width=0.5mm] (5,-4)--(10,-4);

\node (X) at (2,-3.5) {exp decay}; 

\end{tikzpicture}
\end{center} 
\caption{The phase diagrams of the loop $\mathrm{O}(1)$, single random current, and random-cluster measures on the $d$-regular wired tree $\T^d$ coincide. The free measures on $\mathtt{C}_{n}^{d}$, the $d$-regular tree where every edge is substituted by a cycle, have different phase transitions. Finally, the free loop $\mathrm{O}(1)$ model on the monster $\M$ (constructed in the proof of Theorem \ref{thm:non_uniqueness}) has a non-unique phase transition. This is to be contrasted with the corresponding table for the hypercubic and hexagonal lattices in \cite[Figure 1]{hansen2023uniform}.  \label{fig:Phase diagram}}

\end{figure}

\subsection{The graphical representations of Ising} \label{sec:definitions}
We define the random-cluster and random current model as in \cite{klausen2021monotonicity,klausen2023random} through the couplings to the loop $\mathrm{O}(1)$ model. 
Given a finite graph $G = (V,E),$ an \emph{even subgraph} $(V,F)$ of $G$ is a spanning subgraph where each $ v \in V $  is incident to an even number of edges in $F$.  We let $ \Omega_\emptyset(G)$ denote the set of even subgraphs of $G$.
The loop $\mathrm{O}(1)$ model $\ell_{x,G}$ is a natural probability measure on $ \Omega_\emptyset(G) $:
\begin{align}
\ell_{x,G}[\eta] = \frac{1}{Z_G} x^{\abs{\eta}}, \text{     for each } \eta  \in \Omega_\emptyset(G)
\end{align}
with $Z_G = \sum_{\eta \in \Omega_\emptyset(G)} x^{\abs{\eta}} $. 
 Here $\abs{\eta}$ denotes the number of (open) edges in $\eta$ and $x = \tanh(\beta)\in (0,1)$ as in \cite{Lis}, where $\beta$ is the inverse temperature. For $G$ a graph with boundary, we denote $\ell^1_{x,G}=\ell_{x,G/\sim},$ where $\sim$ identifies the boundary vertices of $G$. We refer to $\ell^1_{x,G}$ as the \textit{wired} loop O($1$) model.

We denote Bernoulli edge percolation with parameter $x \in [0,1]$ by $\mathbb{P}_x$ and define the (traced, sourceless) \textit{single random current} at parameter $x$ and boundary condition $\xi\in \{0,1\}$ as
\begin{align}
\Prbcur^{\xi}_x = \ell^{\xi}_x \cup \mathbb{P}^{\xi}_{1- \sqrt{1-x^2}},
\end{align}
where $\ell^{0}_x=\ell_x$. This definition of the model is equivalent to the standard definition as explained in \cite[Remark 3.4]{duminil2016random}, see also \cite{lupu2016note}.
 
 Similarly, we define the random-cluster model via
 \begin{align} \label{eq: random-cluster model}
 \phi^{\xi}_x = \ell^{\xi}_x \cup \mathbb{P}^{\xi}_x,
 \end{align}
 which is equivalent to the standard definition of the model by the results \cite{evertz2002new,grimmett2007random}.

The random-cluster model satisfies several useful monotonicity properties which are not enjoyed by the loop O($1$) and random current models \cite{klausen2021monotonicity}. We endow $\{0,1\}^E$ with the pointwise partial order $\preceq$, and say that an event $A$ is increasing if $\omega\in A$ and $\omega\preceq \omega'$
implies that $\omega'\in A$. The following monotonicity properties will be of use in this paper (see \cite[Theorem 1.6]{DC17}):
\begin{enumerate}
    \item  The FKG inequality: $\phi_{x,G}[A\cap B]\geq \phi_{x,G}[A]\phi_{x,G}[B]$ for $A$ and $B$ increasing.
    \item Monotonicity in boundary conditions: $\phi_{x,G}^1[A]\geq \phi_{x,G}^0[A]$ for $A$ increasing.
    \item Stochastic monotonicity: $\phi_{x_2,G}[A]\geq \phi_{x_1,G}[A]$, whenever $x_1<x_2$ and $A$ is increasing. We write $\phi_{x_1,G}\preceq \phi_{x_2,G}$.
\end{enumerate}

The last property is equivalent to the existence of an increasing coupling - that is, a probability measure $\mu$ with marginals $\omega_1 \sim \phi_{x_1}$ and $\omega_2 \sim \phi_{x_2}$ such that $\omega_1 \preceq \omega_2$ almost surely (i.e.\ $\omega_1$ is a subgraph of $\omega_2$).

We refer to the lecture notes of Duminil-Copin \cite{DC17} for an overall introduction to the Ising model and its graphical representations. Furthermore, we stick with the parametrisation in terms of the loop O($1$) parameter $x \in \lbrack 0, 1\rbrack$ throughout this paper and refer to \cite[Table 1]{hansen2023uniform} for an overview of the standard parametrisations.

\subsection{Graphical representations and uniform even subgraphs.}
In the following, we will consider the uniform even subgraph $\UEG$, which not only serves as the limiting case of the loop $O(1)$ model $\ell_x$ for $x=1$, but also yields (perhaps surprising) connections between the different graphical representations.
For a finite graph $G,$ the uniform even subgraph $\UEG_{G}$ is a uniform element of $\Omega_{\emptyset}(G)$, the set of even subgraphs of $G$.
In \cite{hansen2023uniform}, an abstract view of the uniform even subgraph was taken as the Haar measure on the group of even graphs\footnote{With the group operation given by pointwise addition mod 2 in the space $\{0,1\}^E$ - or, equivalently, taking symmetric differences of edge sets.}, and its percolative properties were studied.  Before that, the uniform even subgraph and its infinite volume measures were studied in detail by Angel, Ray and Spinka in \cite{angel2021uniform}, where the free $\UEG^{0}$ and wired uniform even subgraphs $\UEG^1$ were introduced and it was shown that they coincide on one-ended graphs \cite[Lemma 3.9]{angel2021uniform}. 
In this article, we are concerned with tree-like graphs (as opposed to, say, $\mathbb{Z}^d$), which in general have infinitely many ends. Hence, the distinction between the free and wired measures, $\UEG^{0}$ and $\UEG^1,$ plays a bigger role than in \cite{angel2021uniform,hansen2023uniform}.

For an infinite graph $\mathbb{G}$ the wired uniform even subgraph $\UEG_{\mathbb{G}}^1$ can be defined as the Haar measure (normalised to probability) on $\Omega_{\emptyset},$ the group of all even subgraphs of $\mathbb{G}$. In particular, it pushes forward to Haar measures under group homomorphisms and as a consequence, its marginals are also Haar measures on their supports.

The set of all finite even graphs $\Omega_\emptyset^{< \infty}(\mathbb{G}) = \{ \eta \in \Omega_{\emptyset}(\mathbb{G}) \mid \abs{\eta} < \infty \}$ is a subgroup of $\Omega_{\emptyset}$. Its closure
$$
\Omega^{0}(\mathbb{G}) = \overline{\Omega_\emptyset^{< \infty}(\mathbb{G})}
$$
is  a (compact) subgroup of $\Omega_{\emptyset}(\mathbb{G})$, and the free uniform even subgraph  $\UEG_{\mathbb{G}}^0$ is defined to be the Haar measure on that group. For more details on the construction of the free and wired uniform even subgraphs for infinite graphs and a characterization of the Gibbs of the uniform even graph on trees see \cite[Section 3.2]{hansen2023uniform}.

In \cite[Theorem 3.5]{grimmett2007random} it was realised that the loop $\mathrm{O}(1)$ model arises as the uniform even subgraph of the random-cluster model. Similarly, loop $\mathrm{O}(1)$ model  arises as the uniform even subgraph of the double random current which appeared implicitly in \cite[pf. of Theorem 3.2]{lis2017planar} and was proven in \cite[Theorem 4.1]{klausen2021monotonicity}. 

Thus, on any graph, the loop $O(1)$ measure can be written as follows:
\begin{align}\label{eq:loop-cluster-coupling}
\ell_{\mathbb{G}}^{\xi}[\,\cdot\,] = \phi^{\xi}_{\mathbb{G}}\Big[\text{UEG}^{\xi}_\omega[\,\cdot\,]\Big],
\end{align}
where $\xi = 0$ in the free case and $\xi = 1$ in the wired case\footnote{For finite graphs we write $\ell$, omitting the boundary condition.}, and $\phi^{\xi}_{\mathbb{G}}$ is defined as a thermodynamic limit (see e.g. \cite[Sec. 2.1.3.]{hansen2023uniform}) when $\mathbb{G}$ is infinite. In infinite volume, this may be taken as the definition of the loop O($1$) model (cf. \cite[(4)]{hansen2023uniform}). 
This can then be used to define the single and double random current,\begin{align} \label{Current_couple}
\Prbcur_x^\xi = \ell_x^\xi \cup \mathbb{P}_{1- \sqrt{1-x^2}}
\quad\text{and}\quad
\Prbcur_x^\xi \cup \Prbcur_x^\xi =\ell_x^\xi \cup \ell_x^\xi \cup \mathbb{P}_{x^2}. 
\end{align}

\subsection{Percolation regimes}
Since the graphs we will work on are not vertex-transitive, we will use the following definition of percolation:
We say that a percolation measure $\mu_{x,\mathbb{G}}$ on an infinite graph $\mathbb{G}$ \emph{percolates} if $\mu_{x,\mathbb{G}}[\mathcal{C}_\infty] > 0$ (recall that $\mathcal{C}_\infty$ denotes the event that there exists an infinite cluster) and we define the \emph{percolation regime}
\begin{align}\label{eq:percolation_regime}  
\mathcal{P}(\mu_{x,\mathbb{G}}) = \{x \in (0,1) \mid 
\mu_{x,\mathbb{G}}[\mathcal{C}_\infty] > 0 \}.
\end{align}
We say that the phase transition on $\mathbb{G}$ is unique if both  $\mathcal{P}(\mu_{x,\mathbb{G}})$ and $(0,1)\setminus\mathcal{P}(\mu_{x,\mathbb{G}})$ are connected. In that case, we define the critical parameter $x_c(\mu_{x,\mathbb{G}}) = \inf \mathcal{P}(\mu_{x,\mathbb{G}})$. 
By stochastic monotonicity, the phase transition of the random-cluster model $\phi$ is unique on any graph.

\subsubsection{Bernoulli percolation on a tree.}
We denote by $\mathbb{T}^d$ the $d$-regular tree and by $\mathbb{T}^d_n$ the ball of size $n$ for the graph distance on $\mathbb{T}^d$.
Observing that the cluster of the origin can be described in terms of a Galton-Watson process (and with the observation that vertex-transitivity implies that $\P_{p,\mathbb{T}^d}[0 \cc \infty] > 0$ if and only if $\P_{p,\mathbb{T}^d}[\mathcal{C}_\infty]=1$), one sees that the critical parameter for Bernoulli percolation on the $d$-regular tree is 
\begin{equation}
p_c(\mathbb{P}_{\mathbb{T}^d})=\frac{1}{d-1}. \label{thm:Bernoulli_on_tree}
\end{equation}

\subsubsection{Percolation versus long-range order} 
For readers more familiar with models on lattices, a brief word of caution might be in order: One might wonder why models that essentially share correlation functions nonetheless have different critical parameters for percolation.

In particular, we have the following agreement of two-point functions (see e.g. \cite[Corollary 1.4, Lemma 4.3]{DC17}): 
\begin{equation} \label{Graph rep}
\phi^0_{x,G}[v\cc w]^2= \langle \sigma_v \sigma_w \rangle_{G,\beta}^2  = \Prbcur^0_{x,G}\cup \Prbcur^0_{x,G}[v\cc w]
\end{equation}
for all finite graphs $G$ and vertices $v,w$ (here $\langle \sigma_v \sigma_w \rangle_{G,\beta}$ is the Ising correlation function and $\beta = \arctanh(x)$). However, this is not an obstruction to percolation setting in at different values of $x$ because percolation does not, in general, imply anything for the two-point function. One is tempted to write that the bound $\phi^0[v\cc w]\geq \phi^0[v\cc \infty]\phi^0[w\cc\infty]$ follows from the FKG inequality, but this only holds if the infinite cluster is unique, which, as may be checked, is never true for trees. 

\section{Non-uniqueness of percolation}
In this section we tackle the following question. 

\begin{question} \label{question:mon}
Let $G\subseteq G'$ be two infinite graphs and suppose that the uniform even subgraph of $G$ almost surely percolates. Does the uniform even subgraph of $G'$ almost surely percolate? 
\end{question}

The following counterexample answers the question negatively and at the same time constructs a class of graphs for which $\mathcal{P}(\ell_{x})$ has multiple connected components.

We use that the loop $O(1)$ model factorises on graphs lacking certain cycles:
A cycle denotes a path of vertices $v_1,v_2,...,v_n$ such that $v_1=v_n$. We say that a cycle is \emph{simple} if $v_j\neq v_k$ for distinct $1\leq j,k<n$.

\begin{definition}
For a graph $G$ and two subgraphs $G_1,G_2\subseteq G$, we say that $(G_1,G_2)$ is a \textbf{cut-point factorisation} of $G$ if  $E(G) = E(G_1) \dot \cup E(G_2)$ and it holds that there is no simple cycle in $G$ which contains edges from both $E(G_1)$ and $E(G_2)$.
\end{definition}

\begin{definition}
We say that a graph-indexed family of percolation measures $\nu_{G}$ \textbf{cut-point factorises} if $\nu_{G}=\nu_{G_1}\otimes \nu_{G_2}$ whenever $(G_1,G_2)$ is a cut-point factorisation of $G$.
\end{definition}
The following lemma will be useful. 
\begin{lemma}[Cut point factorisation] \label{loop_factor} 
Each of the measures $\ell_x, \phi_x, \Prbcur_x, \Prbcur_x \cup\Prbcur_x$ cut-point factorises.
\end{lemma} 
\begin{proof}
We first prove the statement for the loop $\mathrm{O}(1)$ model.  For any even subgraph $\eta$ of $G,$ its restrictions $\eta_1, \eta_2$ to the subgraphs $G_1$ and $G_2$ are even graphs since $(G_1,G_2)$ is a cut-point factorisation of $G$.  Thus, writing $\Omega_i=\{0,1\}^{E_i},$ one may simply rewrite
$$ \ell_{x,G}[\eta]= \frac{x^{\abs{\eta}}}{\sum_{\eta' \in \Omega} 1_{\partial \eta' = \emptyset} x^{\abs{\eta'}}} = \frac{x^{\abs{\eta_1}} x^{\abs{\eta_2}}}{\sum_{\eta'_i \in \Omega_i} 1_{\partial \eta'_1 = \emptyset}1_{\partial \eta'_2 = \emptyset} x^{\abs{\eta'_1}} x^{\abs{\eta'_2}}} = \ell_{x,G_1}[\eta_1]\ell_{x,G_2}[\eta_2].
$$
For the other models, the statement follows from the couplings of $\Prbcur, \Prbcur \cup \Prbcur$ and $\phi$ to the loop $\mathrm{O}(1)$ model (which we used to define said models in Section \ref{sec:definitions}).
Indeed, note that for any two product measures, $\mu_{x,G} = \mu_{x,G_1} \otimes \mu_{x,G_2}$ and $\nu_{x,G} = \nu_{x,G_1} \otimes \nu_{x,G_2}$ that
$$
\mu_{x,G} \cup \nu_{x,G}  =  ( \mu_{x,G_1} \cup \nu_{x,G_1}) \otimes (\mu_{x,G_2} \cup \nu_{x,G_2}),
$$
and Bernoulli percolation clearly cut-point factorises.
\end{proof}
In \Cref{loop_factor} we proved the cut-point factorisation property for finite graphs. The infinite volume measures also have the cut-point factorisation property as they are either limits or (sprinkled) uniform even subgraphs of a cut-point factorising measure. 

\begin{lemma}[Non-monotonicity of loop $O(1)$ two-point function] \label{loop_non_monotonicity} 
There exist parameter values $0 < x_1 < x_2 < 1$, and a finite graph $G^{\diamond}$ with vertices $a,b$ such that
$$
\ell_{x_2, G^{\diamond}}[a \cc b] < \frac{1}{4} < \ell_{x_1, G^{\diamond}}[a \cc b]. 
$$
\end{lemma}
\begin{proof} 
Consider the graph $G^{\diamond}$, described in the right-most column of Figure \ref{counter} (which was previously used as an example in \cite{klausen2021monotonicity}). 
In this graph, the probability that the two marked vertices are connected is 
$$
\ell_{x, G^{\diamond}}[a \cc b] = \frac{x^{2m}+x^{2m+2n}}{1+x^{2n}+x^{2m} + 4x^{n+m} + x^{2m+2n}}. 
$$
Setting $n=12,m=2,x_1=0.85$, and $x_2 = 0.965$ yields
$$
\ell_{x_2, G^{\diamond}}[a \cc b]  \leq 0.245 < \frac{1}{4}, \text{ and } \ell_{x_1, G^{\diamond}}[a \cc b]    \geq 0.27 > \frac{1}{4}. 
$$
We have also plotted the graph of $x \mapsto \ell_{x, G^{\diamond}}[a \cc b] $ in Figure \ref{fig:graph_example}.
\end{proof}

	\begin{figure}
		{\caption{The graph $G^{\diamond}$ (pictured to the right) along with its eight even subgraphs (including $G^{\diamond}$ itself).  We let the outer paths be $n$ edges long and the inner paths be $m$ edges long. The nodes $a$ and $b$ are marked with dots. We list the number of edges of each subgraph, the corresponding weights and whether $a$ and $b$ are connected in the subgraph.  (Sketch and text partially revised from \cite{klausen2021monotonicity}.) \label{counter} }}
		
		\begin{center}
		{ 	\begin{tabular}{c |c | c | c | c | c | c | c| c  } 
				\hline
			  Subgraph &	\begin{tikzpicture}[scale = 0.5]
			\draw (1,1) node[fill, circle, scale=0.25] {.};
			\draw (1,-1) node[fill, circle, scale=0.25] {.};
			\draw (1,-2) ; 
			\end{tikzpicture}
			& 
			\begin{tikzpicture}[scale = 0.5]
			\draw  (0,0) -- (0,2) --(2,2)  -- (2,0)  (0,0) -- (0, - 2) --(2, - 2)  -- (2,0)  ;
			\draw (1,1) node[fill, circle, scale=0.25] {.};
			\draw (1,-1) node[fill, circle, scale=0.25] {.};
			\end{tikzpicture}
			&
			\begin{tikzpicture}[scale = 0.5]
			\draw (0,0) --(1,1) -- (2,0)  (0,0) -- (1,-1) -- (2,0)  ;
			\draw (1,1) node[fill, circle, scale=0.25] {.};
			\draw (1,-1) node[fill, circle, scale=0.25] {.};
			\draw (1,-2) ;
			\end{tikzpicture}
			&
			\begin{tikzpicture}[scale = 0.5]
			\draw   (0,0) -- (1,-1) -- (2,0)  (0,0) -- (0,2) --(2,2)  -- (2,0)    ;
			\draw (1,1) node[fill, circle, scale=0.25] {.};
			\draw (1,-1) node[fill, circle, scale=0.25] {.};
			\draw (1,-2) ;
			\end{tikzpicture}
			&
			\begin{tikzpicture}[scale = 0.5]
			\draw   (0,0) -- (1,-1) -- (2,0)   (0,0) -- (0, - 2) --(2, - 2)  -- (2,0)  ;
			\draw (1,1) node[fill, circle, scale=0.25] {.  };
			\draw (1,-1) node[fill, circle, scale=0.25] {. };
			\end{tikzpicture}
			&
			\begin{tikzpicture}[scale = 0.5]
			\draw (0,0) --(1,1) -- (2,0)   (0,0) -- (0, - 2) --(2, - 2)  -- (2,0)  ;
			\draw (1,1) node[fill, circle, scale=0.25] {.};
			\draw (1,-1) node[fill, circle, scale=0.25] {.};
			\end{tikzpicture}
			&
			\begin{tikzpicture}[scale = 0.5]
			\draw (0,0) --(1,1) -- (2,0)    (0,0) -- (0,2) --(2,2)  -- (2,0)   ;
			\draw (1,1) node[fill, circle, scale=0.25] {.};
			\draw (1,-1) node[fill, circle, scale=0.25] {.};
			\draw (1,-2) ;
			\end{tikzpicture}
			&
			\begin{tikzpicture}[scale = 0.5]
			\draw (0,0) --(1,1) -- (2,0)  (0,0) -- (1,-1) -- (2,0)  (0,0) -- (0,2) --(2,2)  -- (2,0)  (0,0) -- (0, - 2) --(2, - 2)  -- (2,0)  ;
			\draw (1,1) node[fill, circle, scale=0.25] {.};
			\draw (1,-1) node[fill, circle, scale=0.25] {.};
			\draw (1,1.3) node {\tiny{$m$}};
			\draw (1,-1.3) node {\tiny{$m$}};
			\draw (1,2.2) node {\tiny{$n$}};
			\draw (1,-1.8) node {\tiny{$n$}};
			\end{tikzpicture} \\
			\hline 
		 	Edges & 0 & $2n$ &  $2m$ & $n+m$ & $n+m$ & $n+m$ & $n+m$ &  $2m +2n$ \\
			\hline 
			Weight & 1 & $x^{2n}$ &  $x^{2m}$ & $x^{n+m}$ & $x^{n+m}$ & $x^{n+m}$ & $x^{n+m}$ &  $x^{2m +2n}$ \\
				\hline 
			$\{a \cc b \} $ & $\times$  & $\times $ &  \checkmark &  $\times $ &  $\times $ &  $\times $ &  $\times $ & \checkmark\\
				\hline 
			\end{tabular}}
		\end{center}
	\end{figure}

\begin{figure}
    \centering
    \includegraphics[scale = 0.5]{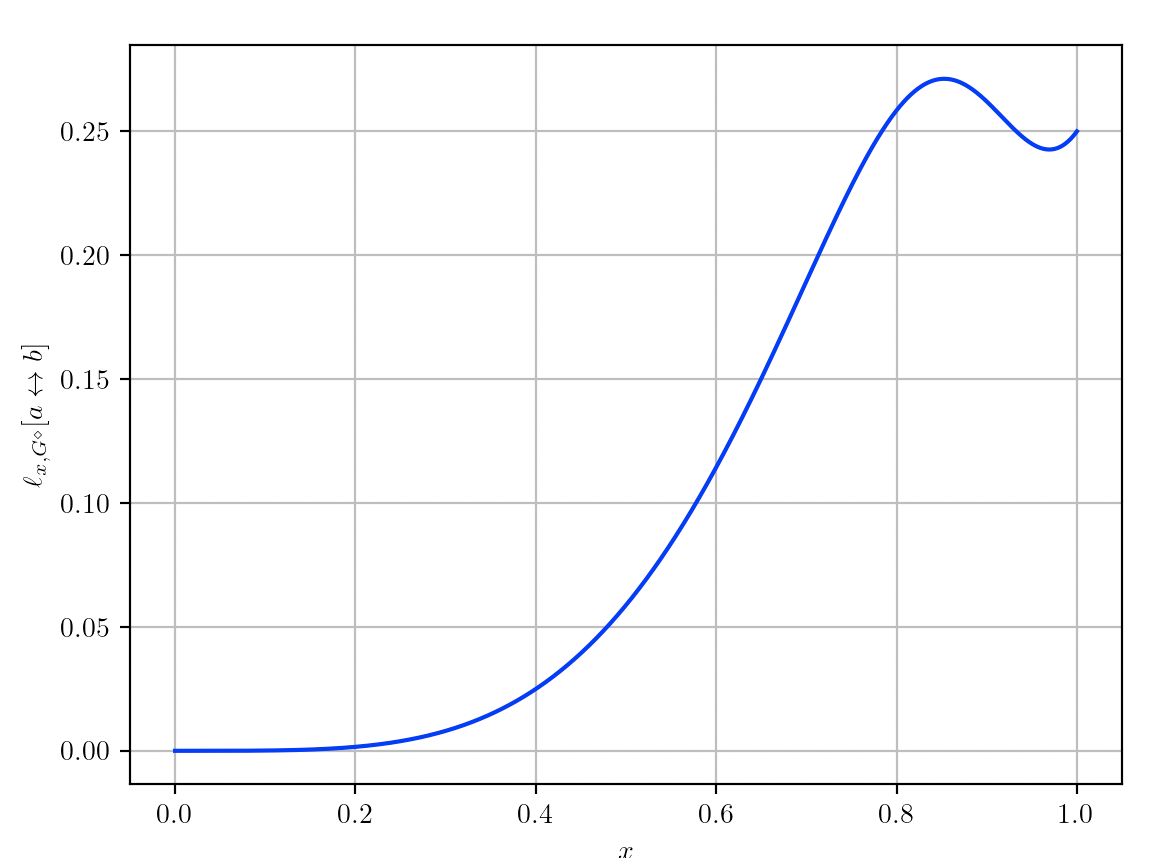}
    \caption{Graph of the connection probability for the loop $\mathrm{O}(1)$ model on the graph $G^{\diamond}$, described in Figure \ref{counter}, for $n=12$ and $m=2$. See also \cite[Figure 2.3]{klausen2021monotonicity} for similar figures.}
    \label{fig:graph_example}
\end{figure}

\begin{figure}
  \centering
  \includegraphics[scale =1]{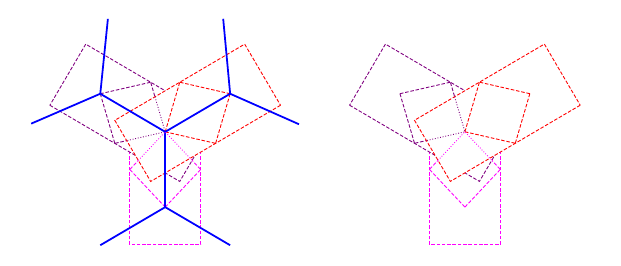}
    \caption{ An illustration of a part of the graph $\M$ built from the graphs $G^{\diamond}$ (see \Cref{counter}) when $d=2$. To the left with the $d$-regular tree overlaid. To the right, the geometry of $\M$ at a single vertex of the initial tree. \label{Monster} }
\end{figure}

With this in hand, we can construct our counterexample: 
\nonuniqueness*
We note that the proof of the theorem also implies that the percolation regime $\mathcal{P}(\ell^0_{x, \mathbb{G}}) \subset [0,1]$ is not connected. 
\begin{proof}
We construct $\mathbb{M}$ from the $d$-regular tree $\mathbb{T}^d$ with root $0$ for an appropriate choice of $d$. Consider the natural orientation of the tree where edges are oriented away from the root and replace every such oriented edge $e=(v,w)$ by a copy of the graph $G^{\diamond}$ where $a$ is identified with $v$ and $b$ is identified with $w$. More formally, let $\mathbb{M}=(\coprod_{e \in E (\mathbb{T}^d)} G_e^{\diamond})/{\sim},$ where $\sim$ is the equivalence relation such that $ a_e\sim b_{e'}$ whenever the source of $e$ is equal to the sink of $e'$. See \Cref{Monster} for an illustration. 
 Now, since the macroscopic structure of $\mathbb{M}$ is that of a tree, Lemma \ref{loop_factor} applies, and it holds that 
\begin{equation} \label{Monster factor}
\ell^0_{x,\mathbb{M}} = \otimes_{e \in E(\mathbb{T}^d) } \ell^0_{x,G^{\diamond}_e}. 
\end{equation}
As a consequence, analysing percolation of $\ell_{x,\mathbb{M}}$ just boils down to Bernoulli percolation on $\mathbb{T}^d$: For given $\eta\in \Omega^0_{\emptyset}(\mathbb{M})$ and $e\in E(\mathbb{T}^d),$ we define $m_e=1$ if  $a \overset{\eta}\cc b$ in $G^{\diamond}_e$ and $m_e = 0$ otherwise.
In other words, $m = (m_e)_{e\in E(\mathbb{T}^d)}$ maps a percolation configuration in $\mathbb{M}$ onto one in $\mathbb{T}^d$ such that $0\, \smash{\overset{m(\eta)}{\longleftrightarrow}}\, \infty $ if and only if $0\, \smash{\overset{\eta}{\longleftrightarrow}}\, \infty $. 
Furthermore, by \eqref{Monster factor}, the image measure of $\ell_{x,\mathbb{M}}$ is Bernoulli percolation: $m(\ell^0_{x,\mathbb{M}}) =  \mathbb{P}_{f(x),\mathbb{T}^d}$ where $f(x) = \ell_{x,G^{\diamond}}[a \cc b]$.

As a consequence, it holds that 
$\ell_{x,M}^0[\mathcal{C}_\infty ] = \mathbb{P}_{f(x),\mathbb{T}^d}[\mathcal{C}_\infty] $. 
By \eqref{thm:Bernoulli_on_tree}, we know that $\mathbb{P}_{f(x),\mathbb{T}^d}[\mathcal{C}_\infty]>0$ if and only if $f(x) > \frac{1}{d+1}$. 
Now, by Lemma \ref{loop_non_monotonicity} there exist $x_1<x_2$ and $d$ such that $f(x_2) < \frac{1}{d+1} < f(x_1)$, which proves the desired.
\end{proof}
As a result, we obtain a negative answer to Question \ref{question:mon}. 
\counterexample*
\begin{proof}[Proof of Corollary~\ref{cor:counterexample} (in the free case)]
First, recall that the loop $\mathrm{O}(1)$ model can be obtained by sampling a uniform even subgraph from a random-cluster configuration, see \eqref{eq:loop-cluster-coupling}.
For the random-cluster model $\phi^0_{x,\mathbb{M}}$ we consider the increasing coupling in $x$. For $x_1 < x_2$ as above, the uniform even subgraph of $\phi^0_{x_1,\mathbb{M}}$ almost surely percolates and the uniform even subgraph of $\phi^0_{x_2,\mathbb{M}}$ almost surely does not.
Since the coupling is increasing, there must exist at least one pair $\omega_1 \preceq \omega_2$  where the uniform even subgraph of $\omega_1$ percolates while that of $\omega_2$ does not. 
\end{proof}

\subsection{Corollary \ref{cor:counterexample} in the wired case.}
In this section, we will construct a supergraph\footnote{One may note that our trick is basically to pick a one-ended supergraph of a two-ended graph. One may check that the wired UEG of a graph with multiple ends percolates (see also \cite[Eq. (14)]{hansen2023uniform}). 
Therefore, our construction does not have the flavour of an optimal solution to the problem.} of $\mathbb{Z}$ for which the wired UEG does not percolate.

Define $\mathbb{G}=(\mathbb{Z},\hat{\mathbb{E}}),$ where $(n,m)\in \hat{\mathbb{E}}$ if either $|n-m|=1$ or $n=-m$. Let $\hat{\mathbb{E}}_{\text{arc}}:=\hat{\mathbb{E}}\setminus \mathbb{E}(\mathbb{Z})$ and $\mathbb{G}_{\text{arc}}=(\mathbb{Z},\hat{\mathbb{E}}_{\text{arc}})$.

\begin{lemma} \label{exercise}
The marginal $\UEG^1_{\mathbb{G}}|_{\hat{\mathbb{E}}_{\text{arc}}}=\mathbb{P}_{\frac{1}{2},\mathbb{G}_{\operatorname{arc}}}.$
\end{lemma}
\begin{proof} This follows from the fact that the restriction group homomorphism $\Even(\mathbb{G})\to \{0,1\}^{\hat{\mathbb{E}}_{\text{arc}}}$ is surjective. To see that the map is surjective, it suffices that all single-edge configurations lie in the image. For a given $e=(-n,n)\in \hat{\mathbb{E}}_{\text{arc}},$ we see that $\id_e$ is the image of the loop containing $[-n,n]\cap \mathbb{Z}$ and $e$.
\end{proof}

\textit{Proof of Corollary \ref{cor:counterexample} (for the wired case).} Let $\mathbb{G}$ be as above. Let $v\in \mathbb{Z}$ and note that on the event that there is an open edge in $\hat{\mathbb{E}}_{\text{arc}}$ outside of $[-|v|,|v|],$ the cluster of $v$ has to be finite. By Lemma \ref{exercise}, with probability 1, infinitely many edges in $\hat{\mathbb{E}}_{\text{arc}}$ are open in $\UEG_{\mathbb{G}}$.  By the previous comment, on this event, all clusters are finite. Hence $\UEG_{\mathbb{G}}$ does not percolate. However, $\mathbb{Z}$ has exactly two even subgraphs, $\eta\equiv 0$ and $\eta\equiv 1$, so $\UEG_\mathbb{Z}^1$ does percolate.  \qed

\subsection{Generalisations and non-uniqueness of random current phase transitions}
We can adapt the earlier construction of a disconnected percolation regime from the loop $O(1)$ model to more general cut-point factorising measures:

\begin{proposition} Let $F$ be a finite graph and  $v,w$ be two vertices.  Let  $\{\mu_{x,F}\}_{x\in [0,1]}$ be a family of cut-point factorising percolation measures. Suppose that there exists a finite graph $F$ such that $x \mapsto \mu_{x,F}[v \cc w]$ is not monotone.  Then, there exists an infinite graph $\mathbb{M}$ such that the percolation regime $\mathcal{P}(\mu_{x,\mathbb{M}})$ is disconnected.
\end{proposition}
\begin{proof}
By assumption, there exist $x_1 < x_2 < x_3$ such that $$\max\{\mu_{x_1,F}[v \cc w], \mu_{x_3,F}[v \cc w] \} <  \mu_{x_2,F}[v \cc w].$$ For $p\in (0,1)$ let $\mathbb{G}_p\sim \mathbb{P}_{p,\mathbb{T}^d}$. Notice that almost surely, $x_c(\mathbb{P}_{x,\mathbb{G}_p})=\frac{1}{p(d-1)}.$ By tuning the parameters $p$ and $d$ appropriately, we can make sure that $$\max\{\mu_{x_1,F}[v \cc w], \mu_{x_3,F}[v \cc w] \} < \frac{1}{p(d-1)} <  \mu_{x_2,F}[v \cc w].$$ Now, construct $\mathbb{M}$ by sampling $\mathbb{G}_p$ and substituting each edge of $\mathbb{G}_p$ by a copy of $F$, gluing in the same way as in our construction in the proof of Theorem \ref{thm:non_uniqueness}.
\end{proof}

The proposition shows that we do not need to "fine-tune" the parameters for the transition points to fit the phase transition of the $d$-regular tree.

As a consequence, we can see that the single random current also admits a disconnected percolation regime:
Recall that \cite[Figure 2.3]{klausen2021monotonicity} gives an example of a graph for which connection probabilities are not monotone, hence, by the above proposition, it follows that we get:

\begin{corollary} \label{cor:random_current_non_monotone}
There exists an infinite graph $\mathbb{G}$ such that $\mathcal{P}(\Prbcur_{G})$ is not connected.
\end{corollary}

\begin{remark} We note that since single site connection probabilities are monotone for $\phi_x$ and $\Prbcur_x \cup \Prbcur_x,$ the counterexamples do not work for these models.
\end{remark}

Another model of statistical mechanics which cut-point factorises is the arboreal gas model. The simplest definition of this model in finite volume is as Bernoulli percolation conditioned to be a forest (i.e.\ conditioned not to contain any cycles):
$$
\nu_{\beta,G}[\omega]=\frac{1}{Z_{G,\beta}} \beta^{\abs{\omega}}\id_{\Even(\omega)=\{0\} },
$$
where $0$ denotes the empty graph, which is even.
It is immediate that if $(G_1,G_2)$ is a cut-point factorisation of $G$ and $F_1,F_2$ are subforests of $G_1$ and $G_2$ respectively, then $F_1\cup F_2$ is a subforest of $G$. Thus, $\nu_{\beta,G}$ cut-point factorises. At present, the connection probabilities of the model are conjectured to be monotone in $\beta$ on all finite graphs $G$ \cite[p.2]{halberstam2023uniqueness}, but supposing that they are not, our construction would go through for this model as well, yielding a graph on which the percolation phase transition of the model is not unique.

\section{Phase transitions of the wired models on the $d$-regular tree coincide}

In this section, we prove the first part of Theorem \ref{thm:main_theorem_equal_different} on trees. The main tool in the proof is the observation that all non-empty even subgraphs in $\mathbb{T}^d$ are infinite. In particular, $\ell^1_{x,\mathbb{T}^d}$ is supported on the empty graph if it does not percolate.

First, we will need a small bookkeeping result. Let $\tilde{\mathbb{T}}^d$ denote the rooted $d$-regular tree. Cycles of $\mathbb{T}^d$ containing an open edge $e$ will consist of two infinite paths in isomorphic copies of $\tilde{\mathbb{T}}^d$. It will therefore be important (and easily demonstrated) that changing the degree of a single vertex does not change the phase transition.
\begin{lemma} \label{Stupid detail}
For any $d \geq 2$, we have that
$
x_c(\phi^1_{\mathbb{T}^d})=x_c(\phi^1_{\tilde{\mathbb{T}}^d}). 
$
\end{lemma}
\begin{proof}
Fix one edge $e$ of $\mathbb{T}^d$ and note that for any $x\in (0,1),$ 
\begin{align*}
 \phi^1_{x,\mathbb{T}^d}[\mathcal{C}_{\infty}]= \phi^1_{x,\mathbb{T}^d}[\mathcal{C}_{\infty}|\; \omega_e=1] \, \phi^1_{x,\mathbb{T}^d}[\omega_e=1]\;+\;\phi^1_{x,\mathbb{T}^d}[\mathcal{C}_{\infty}|\; \omega_e=0]\,  \phi^1_{x,\mathbb{T}^d}[\omega_e=0].
\end{align*}
By tail-triviality (see \cite[Theorem 10.67]{Gri06}), $\phi^1_{x,\mathbb{T}^d}[\mathcal{C}_{\infty}]\in \{0,1\}$ and since $\phi^1_{x,\mathbb{T}^d}[\omega_e=0]\in (0,1),$ we conclude that 
$$
\phi^1_{x,\mathbb{T}^d}[\mathcal{C}_{\infty}|\; \omega_e=0]=\phi^1_{x,\mathbb{T}^d}[\mathcal{C}_{\infty}].
$$
Furthermore, $(\mathbb{V}(\mathbb{T}^d),\mathbb{E}(\mathbb{T}^d)\setminus \{e\})$ has two connected components, both of which are isomorphic to $\tilde{\mathbb{T}}^d$. Permitting ourselves a natural abuse of notation, remark that by the Domain Markov Property \cite[p.8]{DC17}, $\phi^1_{x,\mathbb{T}^d}[\,\cdot\, |\; \omega_e=0]=\phi^1_{x,\tilde{\mathbb{T}}^d}\otimes \phi^1_{x,\tilde{\mathbb{T}}^d}.$ Thus, 
$$
\phi^1_{x,\mathbb{T}^d}[\mathcal{C}_{\infty}] = \phi^1_{x,\mathbb{T}^d}[\mathcal{C}_{\infty}|\; \omega_e=0]=\phi^1_{x,\tilde{\mathbb{T}}^d}\otimes \phi^1_{x,\tilde{\mathbb{T}}^d}[\mathcal{C}_{\infty}] =  1- (1-\phi^1_{x,\tilde{\mathbb{T}}^d}[\mathcal{C}_{\infty}])^2,
$$ which finishes the proof.
\end{proof}
    
\begin{theorem}\label{thm:tree_coincide}
For $x>x_c(\phi^1_{\mathbb{T}^d}),$ there exists $c>0$ such that for any vertex $v$, 
$$
\ell^1_{x,\mathbb{T}^d}[v\cc \infty] \geq c.
$$
In particular, $x_c(\ell^1_{\mathbb{T}^d})=x_c(\phi^1_{\mathbb{T}^d})$.
\end{theorem}
\begin{proof} Consider two neighbouring vertices $v,w$ and define the event $A_{v} = \{v \longleftrightarrow \infty \text{ in } \mathbb{T}^d\setminus(v,w)\}$ and analogously $A_{w}$. 
Define the event $L_{v,w} := A_{v} \cap A_{w} \cap \{(v,w)\text{ open}\}$.
In words, this is the event that that there is a loop which contains the edge $(v,w)$ and goes through the wired boundary at infinity.
Now, suppose $x> x_c(\phi^1_{\mathbb{T}^d})$.
In this case, by the FKG inequality, we obtain a lower bound that a random-cluster configuration $\omega$ satisfies $L_{v,w}$:
\[
\phi^1_{x,\mathbb{T}^d}[L_{v,w}] \geq \phi^1_{x,\mathbb{T}^d}[A_{v}]\phi^1_{x,\mathbb{T}_n^d}[A_{w}]\phi^1_{x,\mathbb{T}^d}[(v,w)\mathrm{\; open}]\geq \phi^1_{x,\tilde{\mathbb{T}}^d}[0\cc \infty]^2 x\geq c > 0, 
\]
where, in the second inequality, we have used monotonicity in boundary conditions,  the fact that the probability of a given edge being open is at least $x$ (which follows from \eqref{eq: random-cluster model}), as well as Lemma \ref{Stupid detail}.

Conditionally on $L_{vw}$, there exist two disjoint infinite paths in $\omega$ starting from $v$ and $w$ respectively.
Let us argue that, for such a configuration $\omega,$ $\UEG_{\omega}[ (v,w)\text{ open}, v \cc \infty]=\frac{1}{2}$. This boils down to two observations:  First, because all components of an even subgraph of $\mathbb{T}^d$ are either trivial or infinite, 
\begin{align}\label{eq:ueg-trivial-infinite}
\text{UEG}^1_{\omega}[ (v,w)\text{ open}, v \cc \infty]=\text{UEG}^1_{\omega}[ (v,w)\text{ open}].
\end{align}
Second, the probability of a given edge, which is part of a loop in $\omega$, being open in $\UEG_{\omega}$ is $\frac{1}{2}$ (see the much more general statement \cite[Lemma 3.5]{hansen2023uniform}).

In conclusion, for $x > x_c(\phi^1_{\mathbb{T}^d})$, 
\[
\ell_{x,\mathbb{T}^d}^1[v\cc \infty ]\geq \phi_{x,\mathbb{T}^d}^1[ \mathrm{UEG}^1_{\omega}[v\cc \infty]] > c/2 > 0.
\]
\end{proof}

To finish the proof of the first statement in \Cref{thm:main_theorem_equal_different}, we make a short aside to discuss the subcritical regime of the random-cluster and random current models on the tree. It is classical that, for $x<x_c(\phi_{\mathbb{T}^d}^1),$ we have that $\phi_{x,\mathbb{T}^d}^1=\mathbb{P}_{x,\mathbb{T}^d}$ (see \cite[Theorem 10.67]{Gri06}). A similar result holds for the double random current:

\begin{lemma} \label{Current_FK_compare}
For $x<x_c(\phi_{\mathbb{T}^d}^1),$ then $
\Prbcur_{x,\mathbb{T}^d}^1\cup \Prbcur_{x,\mathbb{T}^d}^1=\mathbb{P}_{x^2,\mathbb{T}^d}.
$
Moreover, $x_c(\phi_{\mathbb{T}^d}^1)\leq x_c(\Prbcur_{\mathbb{T}^d}^1\cup \Prbcur_{\mathbb{T}^d}^1)$.
\end{lemma}
\begin{proof}
Since $x<x_c(\phi_{\mathbb{T}^d}^1)$ and $\ell^1_{x,\mathbb{T}^d}\preceq \phi^1_{x,\mathbb{T}^d}$, we get that
$
\ell_{x,\mathbb{T}^d}^1[\mathcal{C}_{\infty}]\leq \phi_{x,\mathbb{T}^d}^1[\mathcal{C}_{\infty}]=0.
$
 Therefore, since all non-singleton components of an even subgraph of $\mathbb{T}^d$ are infinite, we have that
$$
\ell^1_{x,\mathbb{T}^d}[\eta\equiv 0|\; \Omega_{\emptyset}(\mathbb{T}^d)\setminus \mathcal{C}_{\infty}]=1.
$$
By \eqref{Current_couple}, this implies that $\Prbcur^1_{x,\mathbb{T}^d}\cup \Prbcur^1_{x,\mathbb{T}^d}=\mathbb{P}_{x^2,\mathbb{T}^d}. $ 
For the second statement, note that $ \mathbb{P}_{x^2,\mathbb{T}^d} \preceq \mathbb{P}_{x,\mathbb{T}^d} = \phi^1_{x,\mathbb{T}^d}$.  
\end{proof}

We can now put \eqref{Current_couple}, Theorem \ref{thm:tree_coincide} and Lemma \ref{Current_FK_compare} together to obtain
$$
x_c(\ell^1_{\mathbb{T}^d})\geq x_c( \Prbcur^1_{\mathbb{T}^d}) \geq x_c( \Prbcur^1_{\mathbb{T}^d} \cup\Prbcur^1_{\mathbb{T}^d})\geq x_c(\phi^1_{\mathbb{T}^d})=x_c(\ell^1_{\mathbb{T}^d}).
$$
Hence, we arrive at the following corollary:

\begin{corollary} \label{cor:double_current_tree}
For $d\geq 2$, 
$
x_c(\ell^{1}_{\mathbb{T}^d})=x_c( \Prbcur^{1}_{\mathbb{T}^d}) = x_c( \Prbcur^{1}_{\mathbb{T}^d} \cup\Prbcur^{1}_{\mathbb{T}^d}).
$
\end{corollary}

\medskip
\subsection{Modifications for $\mathtt{C}_n^d$.}
In the following, we comment on how to adapt the previous proof strategy to yield the analogue of Corollary~\ref{cor:double_current_tree} (resp./ the first part of Theorem \ref{thm:main_theorem_equal_different}) on $\mathtt{C}^d_n$. This requires two ingredients:
 \begin{enumerate}
     \item[a)] For $x<x_c(\phi^{1}_{\mathtt{C}_n^d}),$ all models reduce to explicitly comparable independent models. In particular, we will argue that \[x_c(\Prbcur^1_{x,\mathtt{C}_n^d}\cup \Prbcur^1_{x,\mathtt{C}_n^d})\geq x_c(\phi^1_{x,\mathtt{C}_n^d}).\]
     \item[b)] For $x>x_c(\phi^1_{\mathtt{C}_n^d})$, the loop $\mathrm{O}(1)$ model $\ell^1_{x,\mathtt{C}_n^d}$ percolates.
 \end{enumerate}

\begin{onehalfspace} 
 For $a)$ if $x<x_c(\phi^{1}_{\mathtt{C}_n^d}),$ rather than $\eta\sim\ell^1_{x,\mathtt{C}_n^d}$ being deterministically empty, it includes each simple cycle of $\mathtt{C}_n^d$ independently since the free loops cut-point factorise. Accordingly, $\Prbcur^1_{x,\mathtt{C}_n^d}\cup \Prbcur^1_{x,\mathtt{C}_n^d}$ is a union of two independent cycle measures and a Bernoulli percolation and therefore, it percolates only if it has better connection probabilities in finite volume than $\phi^1_{x,\mathtt{C}_n^d}$. But the finite-volume two-point function of $\phi^1_{x,\mathtt{C}_n^d}$ is always larger than that of $\Prbcur^1_{x,\mathtt{C}_n^d}\cup \Prbcur^1_{x,\mathtt{C}_n^d}$ by \eqref{Graph rep}. Since $x<x_c(\phi^{1}_{\mathtt{C}_n^d}),$ we conclude that $\Prbcur^1_{x,\mathtt{C}_n^d}\cup \Prbcur^1_{x,\mathtt{C}_n^d}$ does not percolate.

 Now, for $b),$ if $x>x_c(\phi^1_{\mathtt{C}_n^d}),$ we want to make an observation that infinite paths can be deduced in $\ell^1_{x,\mathtt{C}_n^d}$ from a local configuration. On a tree, it is true that any open edge is part of an infinite cluster (this is what we used in the proof for $\mathbb{T}^d$, see \eqref{eq:ueg-trivial-infinite}). On $\mathtt{C}_n^d,$ instead, it is true that if $e,e'$ are edges belonging to the same simple cycle, and $\eta\in \Even(\mathtt{C}_n^d)$ with $\eta_e=1$ and $\eta_{e'}=0,$ then $e$ lies on an infinite cluster in $\eta$ (in  which case $e'$ and $e$ are on opposite paths between the glued vertices). By the same argument as previously, conditionally on $e$ being cyclic and lying in an infinite component of $\omega\sim \phi^1_{x,\mathtt{C}_n^d},$ the probability that $\eta_e=1$ and $\eta_{e'}=0$ is at least $\frac{1}{4}$ for $\eta\sim \UEG^1_{\omega}$, which concludes the argument.
\end{onehalfspace}
\section{Explicit computation of critical points}
In this section we explicitly compute the critical points for the free models on $\mathtt{C}_{n}^{d}$, the $d$-regular tree where every edge is replaced by a cycle of length $2n$ (and glued through opposite point of the cycle). 

\begin{figure}
    \centering
    \includegraphics[width=0.6\linewidth]{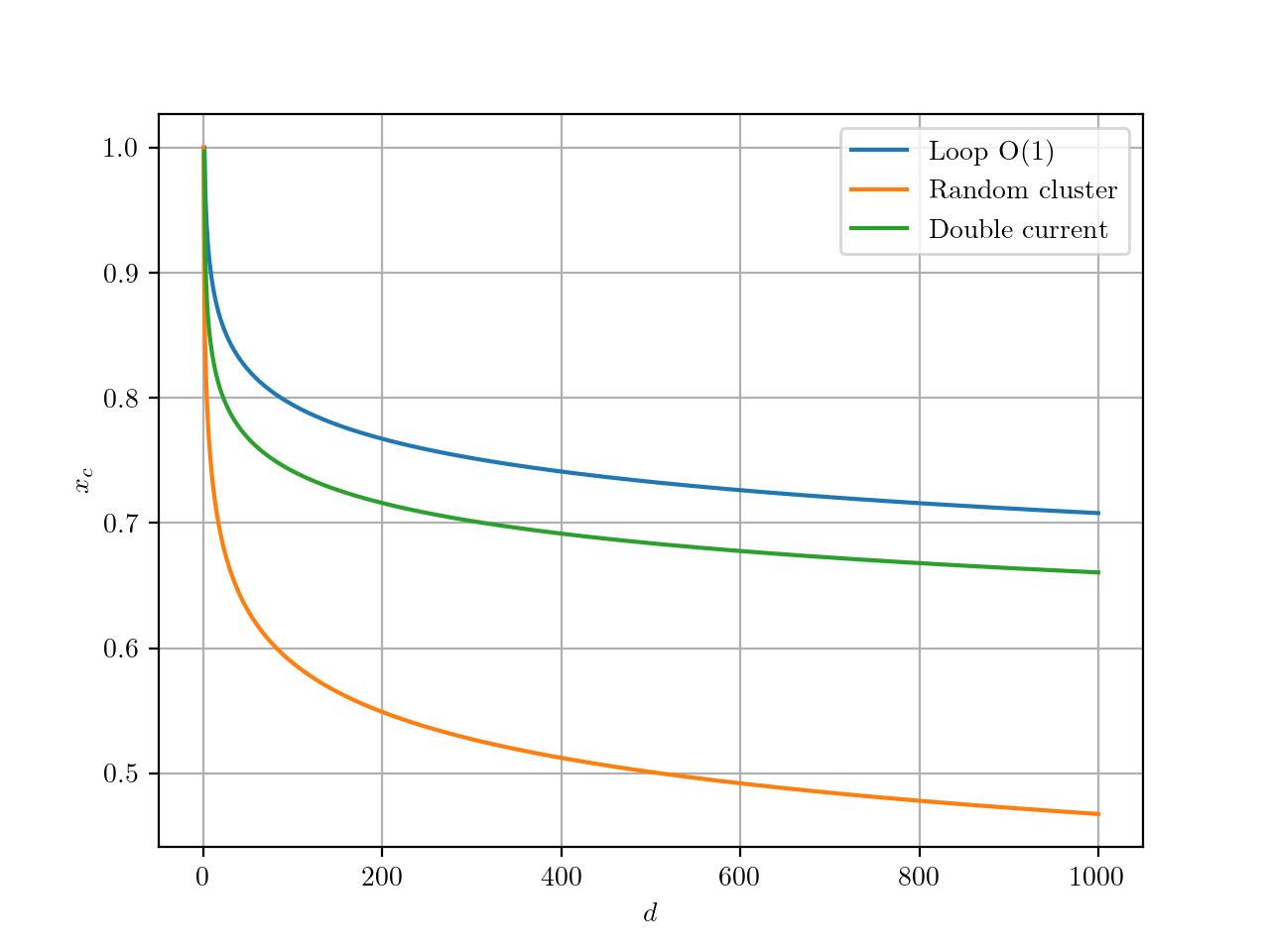}
    \caption{The critical $x_c$ on the graph $\mathtt{C}_{d,n}$ as a function of $d$ for the loop $\mathrm{O}(1)$, and double random current.}
    \label{fig:Critical}
\end{figure}

\begin{proposition} \label{prop:cycle_tree}
  For any $n \geq 1$ and $d \geq 2$, it holds that 
    \begin{align*}
        &x_c(\ell^0_{\mathtt{C}_{n}^{d}}) = (d-2)^{-\frac{1}{2n}} \\
         &x_c( \phi^0_{\mathtt{C}_{n}^{d}} ) = \sqrt[n]{(d-1) - \sqrt{ (d-1)^2 -1}} \\
         &x_c( \Prbcur^0_{\mathtt{C}_{n}^{d}} \cup \Prbcur^0_{\mathtt{C}_{n}^{d}} ) = \sqrt[2n]{(2d-5) - \sqrt{(2d-5)^2 -1} }.
    \end{align*}
    In particular, the three different models have three different phase transitions. 
\end{proposition}
\noindent A graphical presentation of the functions is given in Figure \ref{fig:Critical}.

\begin{proof}
Let $C_{2n}$ denote cycle graph of length $2n$ and let $a$ and $b$ be two antipodal points. 
As in the proof Theorem~\ref{thm:non_uniqueness}, we use cut-point factorisation (Lemma~\ref{loop_factor}) to reduce everything to Bernoulli percolation on $\mathbb{T}^d$ with parameter $\nu_{x, C_{2n}}[a \cc b]$ (with $\nu$ denoting one of the models under consideration). 
The rest follows by direct computation:

\noindent
\textit{Loop $\mathrm{O}(1)$ model:} It holds that 
$
\ell_{x,C_{2n}}[a \cc b]= \frac{x^{2n}}{1+x^{2n}}.
$
Thus,
$$
\frac{1}{d-1} = \frac{x_c^{2n}}{1+x_c^{2n}}
$$
can be solved to obtain $x_c(\ell^0_{\mathtt{C}_{n}^{d}}) = (d-2)^{-\frac{1}{2n}}$.

\medskip\noindent
\textit{Random-cluster model:} 
Now, for $(\eta,\omega)\sim \ell_{x,C_{2n}}\otimes \mathbb{P}_{p,C_{2n}}$, we see that
$$
\ell_{x,C_{2n}}\otimes \mathbb{P}_{p,C_{2n}}[\eta\cup \omega\in (a\cc b)|\;\eta]=\begin{cases}
    1 & \eta\equiv 1\\
    2p^n-p^{2n} & \eta \equiv 0.
\end{cases}
$$
Since the cycle graph has exactly two even subgraphs (the full and the empty graph) we get that
\begin{align}\label{eq:Bersprink}
(\ell_{x, C_{2n}} \cup \Prb_{p, C_{2n}})[a \cc b] =  \frac{x^{2n}+ 2p^n- p^{2n}}{1+x^{2n}}.
\end{align}
Thus, for the random-cluster model (where we choose $p = x$) we obtain
$$
\frac{1}{d-1} = \phi_{x_c,C_{2n}}[a \cc b] = \frac{x_c^{2n}+ 2x_c^n- x_c^{2n}}{1+x_c^{2n}} = \frac{2x_c^n}{1+x_c^{2n}},
$$
which can be solved with the substitution $z =x_c^n$ and
we obtain that
$$
 x_c( \phi^0_{\mathtt{C}_{n}^{d}} ) = \sqrt[n]{(d-1) - \sqrt{ (d-1)^2 -1} }.
 $$
 
\medskip\noindent
\textit{Double random current:} Analogously to \eqref{eq:Bersprink} we obtain
$$
(\ell_{x,C_{2n}} \cup \ell_{x,C_{2n}} \cup \Prb_{p,C_{2n}}) [a \cc b] = \frac{2 x^{2n} + x^{4n} + 2p^{n} - p^{2n} }{(1+x^{2n})^2}. 
$$
Hence, choosing $p = x^2$ for the double random current we obtain the following equation
$$
\frac{1}{d-1} = \frac{4 x^{2n} }{(1+x^{2n})^2}
$$
which can be solved with the substitution $z=x^{2n}$,  giving rise to
$$
x_c( \Prbcur^0_{\mathtt{C}_{n}^{d}} \cup \Prbcur^0_{\mathtt{C}_{n}^{d}} ) = \sqrt[2n]{(2d-5) - \sqrt{ (2d-5)^2 -1} }. 
$$
\end{proof}

The same argument for the single current does not lead to a closed formula, but the following separate argument allows us to conclude \Cref{thm:main_theorem_equal_different}. 

\medskip
\textit{Proof of \Cref{thm:main_theorem_equal_different}}. The first part of the theorem is given by combining \Cref{thm:tree_coincide} and \Cref{cor:double_current_tree}. 
We focus on the single random current:
From \eqref{eq:Bersprink} we sprinkle with $p(x) = 1 - \sqrt{1-x^2}$ to obtain the single current. One may check that for any increasing (differentiable) function $p(x)$ taking values in $[0,1],$ the function $x \mapsto (\ell_x \cup \Prb_{p(x)})[a \cc b] $ is increasing. 
Hence, there exists a unique solution $x_c$ to the equation
$$
\frac{1}{d-1} = \frac{x_c^{2n}+ 2p(x_c)^n- p(x_c)^{2n}}{1+x_c^{2n}}. 
$$
Thus, if $p$ and $q$ are increasing differentiable functions that take values in $[0,1]$  such that $p(x) < q(x)$ for all $x \in [0,1]$ then $x_c( \ell^0_{x,{\mathtt{C}_{n}^{d}}} \cup \P_{p(x),{\mathtt{C}_{n}^{d}}}) > x_c( \ell^0_{x,\mathtt{C}_{n}^{d}} \cup \P_{q(x),\mathtt{C}_{n}^{d}}).$

Using this for the functions $r(x)=0$, $p(x) = 1 - \sqrt{1-x^2}$, and $q(x) = x^2$ together with stochastic domination yields
\begin{align*}
\smash{x_c( \ell_{x,\mathtt{C}_{n}^{d}}^0) > x_c( \underbrace{\ell_{x,\mathtt{C}_{n}^{d}}^0\cup \P_{p(x),\mathtt{C}_{n}^{d}}}_{\Prbcur_{x,\mathtt{C}_{n}^{d}}^0}) > x_c( \ell_{x,\mathtt{C}_{n}^{d}}^0\cup \P_{q(x),\mathtt{C}_{n}^{d}})\geq x_c(\underbrace{\ell_{x,\mathtt{C}_{n}^{d}}^0 \cup \ell_{x,\mathtt{C}_{n}^{d}}^0\cup \P_{q(x),\mathtt{C}_{n}^{d}}}_{\Prbcur_{x,\mathtt{C}_{n}^{d}}^0 \cup\, \Prbcur_{x,\mathtt{C}_{n}^{d}}^0}). }
\end{align*}

\qed

\section{The critical probability for Bernoulli percolation is no obstruction for the UEG}
In the previous section, we considered graphs where the free uniform even subgraph is intimately tied to the behaviour of ordinary Bernoulli percolation on another graph. One might wonder about general links between the behaviour of Bernoulli percolation and that of the UEG. For instance, one might have a suspicion that if a graph $\mathbb{G}$ is easily disconnected in the sense that  the percolation threshold $p_c(\mathbb{P}_{\mathbb{G}})$ is very close to $1$, this lack of connectivity might also impact the UEG. This turns out to be false even for one-ended graphs.

We are going to give two counterexamples:
\begin{itemize}
    \item The first is a construction that can be applied to just about any graph and which admits an easy proof. However, the graphs thus produced are not of inherent interest otherwise.
    \item The second is the infinite cluster of $\phi_{x_c+\varepsilon,\mathbb{Z}^2}$ as $\varepsilon\to 0^+$, which is a more natural object, but for which the proof is more involved.
\end{itemize}
The UEG of the latter model always percolates by \cite[Theorem 1.3]{GMM18}, and it is reasonable to believe that continuity of the phase transition 
implies that breaking even a small fraction $\delta=\delta(\varepsilon)$ of the edges breaks the infinite cluster. However, the proof we give relies on \cite{Super-RSW}, the results of which are (conjecturally) not valid for all planar percolation models with a continuous phase transition (see discussion after \cite[Theorem 7.5]{Super-RSW}). Therefore, the matter is much more subtle than one might expect and restating all necessary prerequisites is beyond the scope of the present paper. As such, we shall settle for referring to suitable places in the literature. However, we find the example to be important in the sense that the graph produced is, in some sense, a much more natural negative resolution to the question.

\subsection{The edge-halving construction}
 For a graph $\mathbb{G}=(\mathbb{V},\mathbb{E}),$ define $\mathbb{G}^{1/2}$ with $\mathbb{V}(\mathbb{G}^{1/2})=\mathbb{V}\cup \mathbb{E}$ and $\mathbb{E}(\mathbb{G}^{1/2})$ consisting of pairs $(v,e)$ with $v\in \mathbb{V}$ and $e$ an edge in $\mathbb{G}$ with $v$ as its one end-point. One may note that $\mathbb{G}^{1/2}$ is bipartite with bi-partition $\mathbb{V}\cup \mathbb{E}$. In pictures, $\mathbb{G}^{1/2}$ is obtained from $\mathbb{G}$ by dividing each edge in two. The point is that doing so does not change the behaviour of the uniform even subgraph at all, while it makes it strictly harder for Bernoulli percolation to percolate.

 We note that since $\mathbb{G}^{1/2}$ is bipartite, a subgraph thereof will have an infinite component if and only if it has a connected component containing infinitely many vertices from $\mathbb{V}.$
\begin{lemma} \label{ueg stability}
For any graph $\mathbb{G}=(\mathbb{V},\mathbb{E})$ there is a group isomorphism $\psi:\Omega_{\emptyset}(\mathbb{G}^{1/2})\to \Omega_{\emptyset}(\mathbb{G})$ such that$$
\smash{\overset{\eta}{v\longleftrightarrow w} \quad \text{if and only if} \quad\overset{\psi(\eta)}{v\longleftrightarrow w}},
$$ for every $\eta\in \Omega_{\emptyset}(\mathbb{G})$ and $v,w\in \mathbb{V}$. 
In particular, $\UEG^1_{\mathbb{G}}$ percolates if and only if $\UEG^1_{\mathbb{G}^{1/2}}$ does.
\end{lemma}
\begin{remark}
We note that the lemma also holds for $\UEG^0$, but omit it from the statement for notational ease. The same proof carries through.
\end{remark}
\begin{proof}
Any $e=(v,w)\in \mathbb{E}$ has degree two in $\mathbb{G}^{1/2}$. Therefore, for any $\eta\in \Omega_{\emptyset}(\mathbb{G}^{1/2}),$ $(v,e)\in \eta$ if and only if $(w,e)\in \eta$. Accordingly,
$$
\psi(\eta)=\{ e\in \mathbb{E}|\; \mathrm{deg}_{\eta}(e)=2\}
$$
defines an even subgraph of $\mathbb{G}$ with the desired connectivity property. One checks that $\psi$ is a group homomorphism and, furthermore, that its inverse is given by
$$
\psi^{-1}(\eta)=\{(v,e)|\; e\in \eta\}. 
$$
Since $\psi$ is a group homomorphism and $\UEG^1_{\mathbb{G}^{1/2}}$ is the Haar measure on $\Even(\mathbb{G}^{1/2})$, it pushes forward to Haar measure on its image under $\psi$, which is $\UEG^1_{\mathbb{G}},$ since $\psi$ is surjective.
\end{proof}

\begin{lemma} \label{pc down}
For any graph $\mathbb{G}=(\mathbb{V},\mathbb{E}),$ it holds that $
p_c(\mathbb{P}_{\mathbb{G}^{1/2}})=\sqrt{p_c(\mathbb{P}_{\mathbb{G}})}. 
$
\end{lemma}
\begin{proof}
The proof proceeds by coupling $\omega_p\sim \mathbb{P}_{p^2,\mathbb{G}}$ and $\omega_{p,1/2}\sim \mathbb{P}_{p,\mathbb{G}^{1/2}}$ for every $p\in [0,1]$ in such a way that $\omega_p\in \{v\leftrightarrow w\}$ if and only if  $\omega_{p,1/2}\in \{v\leftrightarrow w\}$ for every pair $v,w\in \mathbb{V}$. The coupling itself declares that $e=(v,w)\in \omega_p$ if and only if $(v,e)\in \omega_{p,1/2}$ and $(w,e)\in \omega_{p,1/2}$. The process $\omega_{p}$ thus defined is i.i.d. since $\omega_{p,1/2}$ is, and its marginals may be checked to be Bernoulli variables of parameter $p^2$. The desired connectivity property also follows by construction.
\end{proof}

This allows us to prove \Cref{norelation}.
\norelation*
\begin{proof}
Let $\varepsilon>0$ and let $\mathbb{G}_0=\mathbb{Z}^2,$ the uniform even subgraph of which percolates by \cite{GMM18} and for which $p_c(\mathbb{P}_{\Z^2})=\frac{1}{2}$ by Kesten's Theorem \cite{kesten1980critical}. Inductively, define $\mathbb{G}_{j+1}=\mathbb{G}_j^{1/2}$. By Lemma \ref{ueg stability}, we have that 
$\UEG_{\mathbb{G}_j}[0\leftrightarrow \infty]>0$
for every $j$ and by Lemma \ref{pc down}, we have that $p_c(\mathbb{G}_j)=2^{-2^{-j}}$. Picking $j$ sufficiently large proves the desired.
\end{proof}

\subsection{The infinite cluster of the slightly supercritical random-cluster model.}
For our second example, for $x>x_c$ we let $\mathbb{G}^{x}$ denote the infinite cluster of $\phi_{x,\mathbb{Z}^2}$. By \cite[Theorem 3.1]{GMM18}, we have that the uniform even subgraph of $\mathbb{G}^x$ percolates almost surely. Therefore, we obtain a second proof of Theorem \ref{norelation} if we can prove the following:
\begin{proposition}
Almost surely, under the increasing coupling of $\phi_{x,\mathbb{Z}^2}$, we have
$$
\lim_{x\downarrow x_c}p_c(\mathbb{P}_{\mathbb{G}^{x}})=1.
$$
\end{proposition}
\begin{proof} We start by fixing parameters and notation. Fix $\delta\in (0,1)$ and let $\rho>0$ be small enough. For finite $G\subseteq \mathbb{Z}^2,$ let $(\omega_{G}, \xi_{G,\delta})\sim \phi^1_{x_c,G}\otimes \mathbb{P}_{1-\delta,G}$. Furthermore,
let $R_k=[-k,k]\times [-3k,3k]\cap \mathbb{Z}^2$, $\Lambda_k=[-k,k]^2\cap \mathbb{Z}^2$ and let $\mathscr{C}_k$ denote the event that there is a crossing in $R_k$ between its left and right faces.

We claim that if $k=k(\delta,\rho)$ is large enough, then
\begin{equation} \label{Liggett-preansatz}
\phi^1_{x_c,R_k}\otimes \mathbb{P}_{1-\delta,R_k}[\omega_{R_k}\cap \xi_{\delta,R_k}\in \mathscr{C}_k]<\rho.
\end{equation}
Before we indicate how to prove \eqref{Liggett-preansatz}, let us see how it finishes the proof. As any crossing between $\Lambda_k$ and $\Lambda_{3k}$ must cross at least one rotated translate of $R_k$, monotonicity in boundary conditions and a union bound implies that
\begin{equation} \label{Liggett-ansatz}
\phi^1_{x_c,\Lambda_{3k}}\otimes \mathbb{P}_{1-\delta,\Lambda_{3k}}[\omega_{\Lambda_{3k}}\cap \xi_{\delta,\Lambda_{3k}}\in \{\Lambda_{k}\cc \Lambda_{3k}\}]<4\rho.
\end{equation}
It is well-known that if $\rho$ is sufficiently small, an estimate of the form \eqref{Liggett-ansatz} for some $k$ is enough to imply non-percolation by techniques that go back to \cite{liggett1997domination} (see e.g.\ the proof of \cite[Proposition 2.11] {hansen2023uniform}). However, by continuity of the finite volume measures, \eqref{Liggett-ansatz} remains true if $\omega_{\Lambda_{3k}}$ is replaced by $\tilde{\omega}_{\Lambda_{3k}}\sim \phi^1_{x_c+\varepsilon,\Lambda_{3k}}$ for $\varepsilon$ sufficiently small. Upon inspection, this is the same as saying that $p_c(\mathbb{P}_{\mathbb{G}^x})\geq 1-\delta$ almost surely for $x\in (x_c,x_c+\varepsilon),$ which is what we wanted, since, under the increasing coupling, $x\mapsto p_c(\mathbb{P}_{\mathbb{G}^x})$ is almost surely decreasing.

Now, to see that \eqref{Liggett-preansatz} holds provided $k$ is large enough, we refer to \cite[Lemma 5.2]{manolescuboole}. This lemma is stated in the context of Boolean percolation, but as is remarked on in that paper, its proof only relies on the techniques of \cite[Theorem 7.5]{Super-RSW}. Thus, it is also valid for the random-cluster model. Combining \cite[Lemma 5.2]{manolescuboole} with the fact that the four-arm exponent of the random-cluster model is smaller than 2 \cite[Page 11]{ScalingRelations} yields \eqref{Liggett-preansatz}.
\end{proof}
\section*{Acknowledgements} 
FRK acknowledge the Villum Foundation for funding through the QMATH Center of Excellence (Grant No. 10059) and the Villum Young Investigator (Grant No. 25452) programs and the Carlsberg Foundation, grant CF24-0466. UTH acknowledges funding from Swiss SNF. PW was supported by the
European Research Council under the European Union’s Horizon 2020 research and
innovation programme (grant agreement No. 851682 SPINRG).  We thank Boris Kjær and Franco Severo for many discussions and comments.

\bibliographystyle{plain}
\bibliography{bibliography}

\end{document}